\documentclass[12pt,a4paper]{amsart}
\usepackage{graphicx}
\usepackage{latexsym}
\usepackage{amssymb}
\usepackage[all]{xy}
\newtheorem{prop}{Proposition}[section]
\newtheorem{cor}[prop]{Corollary}
\newtheorem{thm}[prop]{Theorem}
\newtheorem{lemma}[prop]{Lemma}
\newtheorem{dfn}[prop]{Definition}
\newtheorem{mex}[prop]{Example}
\newtheorem{rem}[prop]{Remark}

\newenvironment{pf}{{\bf Proof.}}{\hfill$\Box$}

\renewcommand{\P}{\mathcal{P}}
\renewcommand{\[}{\begin{equation}}
\renewcommand{\]}{\end{equation}}

\newcommand{\co}{{\mathrm{co}H}}

\newcommand{\can}{\mathsf{can}}
\newcommand{\id}{\operatorname{\rm id}}
\newcommand{\tcan}{\widetilde{\mathsf{can}}}
\newcommand{\s}{s}
\newcommand{\Achain}{\mathcal{A}}

\newcommand{\Upper}{\mathfrak{u}}
\newcommand{\Proj}{{{\mathbb P}^{N-1}({\mathbb Z/2})}}

\begin{document}
\title[Piecewise principal comodule algebras]
{\large Piecewise principal comodule algebras}

\author{Piotr~M.~Hajac}
\address{Instytut Matematyczny, 
Polska Akademia Nauk,
ul.~\'Sniadeckich 8, Warszawa, 00-956 Poland\\
Katedra Metod Matematycznych Fizyki, 
Uniwersytet Warszawski,
ul. Ho\.za 74, Warszawa, 00-682 Poland} 
\email{http://www.impan.gov.pl/\~{}pmh,
http://www.fuw.edu.pl/\~{}pmh}
\author{Ulrich~Kr\"ahmer}
\address{Instytut Matematyczny, Polska Akademia Nauk, 
ul.~\'Sniadeckich 8, Warszawa, 00-956 Poland}
\email{kraehmer@impan.gov.pl}
\author{Rainer~Matthes}
\address{Katedra Metod Matematycznych Fizyki, 
Uniwersytet Warszawski,
ul. Ho\.za 74, Warszawa, 00-682 Poland}
\email{matthes@fuw.edu.pl}
\author{Bartosz~Zieli\'nski}
\address{Instytut Matematyczny, Polska Akademia Nauk, 
ul.~\'Sniadeckich 8, Warszawa, 00-956 Poland\\
Department of Theoretical Physics II, University of \L{}\'od\'z,
Pomorska 149/153 90-236 \L{}\'od\'z, Poland}
\email{bzielinski@uni.lodz.pl}

\maketitle

\begin{abstract}
A comodule algebra $P$ over a Hopf algebra 
$H$ with bijective antipode
is called principal if the coaction of $H$ is Galois
and $P$ is 
$H$-equivariantly projective
(faithfully flat) over the coaction-invariant
subalgebra~$P^\co$. 
We prove that  principality  is
a piecewise property:  given $N$
comodule-algebra surjections $P \rightarrow P_i$ whose
kernels 
intersect to zero, $P$ is principal if and only if
all $P_i$'s are principal. Furthermore, assuming the principality
of $P$, we show that the lattice these kernels generate
is distributive  if and only if so is the lattice obtained by
intersection with $P^\co$.
Finally, assuming the above distributivity property, 
we obtain a flabby sheaf 
of principal comodule algebras over 
a certain space that is universal for all such $N$-families
of surjections $P \rightarrow P_i$  and
such that the comodule algebra of global sections is
$P$.
\end{abstract}

\section{Introduction}

We are motivated by studying equivariant pullbacks of $C^*$-algebras
exemplified, for instance, by a join construction for compact quantum
groups. Insisting on equivarince to be given by ``free actions" leads
to Galois extensions of $C^*$-algebras by Hopf algebras, and generalising
pullbacks to ``multi-pullbacks" brings up distributive lattices as a
fundamental language. As a by-product of our considerations, we obtain
an equivalence between  the distributive lattices generated by
 $N$ ideals 
intersecting to zero   and flabby sheaves of algebras over a certain
$(N-1)$-projective space.

Comodule algebras  provide a natural
noncommutative geometry generalisation of spaces
equipped with group actions. Less evidently,
principal extensions \cite{bh04} appear to be 
a proper analogue of principal bundles in this
context (see Section~\ref{prelimi} for precise definitions).
Principal extensions can be considered as functors from
the category of finite-dimensional corepresentations
of the Hopf algebra (replacing the structure group) 
to the category of finitely generated projective 
modules over the coaction-invariant subalgebra
(playing the role of the base space). 

The aim of this article is to establish a viable concept
of locality of comodule algebras and  analyse its relationship
with principality. 
The notion of locality we use herein results from
decomposing algebras into
``pieces'', meaning expressing them as 
multiple fibre products
(called multirestricted direct sums in \cite[p.~264]{p-gk99}).
If $X$ is a compact Hausdorff space and 
$X_1,\ldots,X_N$ form a finite closed covering, then 
$C(X)$ can be expressed as such a multiple fibre
product of its quotient $C^*$-algebras $C(X_i)$. 
This leads to a $C^*$-algebraic notion of
a ``covering of  a  quantum space" 
given by a finite family of algebra 
surjections 
$ \pi_i : P \rightarrow P_i$ with 
$\bigcap_i \mathrm{ker}\, \pi_i=0$ (see \cite{bk96,cm00}, cf.~\cite{d-m97}).

Recall that not all properties of group actions
are local in nature: there is a natural example of a locally proper 
action of ${\mathbb R}$ on ${\mathbb R}^2$ that is not proper 
(see \cite[p.298]{p-rs61}, 
cf.~\cite[Example~1.14]{bhms} for more details).
 On the
other hand, a group action is free if and only if it is locally free. 
Therefore, since for compact groups all actions are proper, the principal
(i.e., free and proper) actions of compact groups are local in nature. 
 Our  main result 
(Theorem~\ref{main}) is
a noncommutative analogue of this statement: 
a comodule algebra $P$ which is covered 
 by ``pieces'' $P_i$ is principal if and only if so are the pieces.
In particular, a smash product 
of an $H$-module algebra $B$ with the Hopf algebra $H$ 
(with bijective antipode) is principal, so that gluing together
smash products is a way of constructing principal comodule
algebras.

However, it was pointed out in 
\cite[p.369]{cm00} that the aforementioned coverings can show a certain
incompleteness when going beyond the 
$C^*$-setting. This is related to the fact that 
the lattice of ideals generated by the 
$ \mathrm{ker}\, \pi_i$'s is in
general not distributive. (This problem does not arise
for $C^*$-algebras.) Hence we analyse
a stronger notion of covering that 
includes this nontrivial assumption as part of the 
definition (see Definition~\ref{covdef}).
 If all $P_i$'s are smash products,
 we arrive at a concept of piecewise
trivial comodule algebras. They appear to be a good
noncommutative replacement of locally trivial compact
principal bundles.

Our motivation for going beyond $C^*$-algebras comes from the way
we consider compact principal bundles (the Hausdorff property assumed).
We aim  to use at the same time algebraic techniques of Hopf-Galois
 theory and analytic tools coming with $C^*$-algebras. To this end,
 we look at the total space
of such a bundle in terms of the algebra of functions continuous 
along the base and polynomial along the fibres \cite{bhms}. 
Then the base space algebra is always a $C^*$-algebra, but,
 unless the 
group is finite, the total space algebra is not $C^*$. 

The data of a covering by $N$ pieces can be
equivalently encoded into a flabby sheaf of algebras
over ${{\mathbb P}^{N-1}({\mathbb Z/2})}$. This is 
the 2-element field $(N-1)$-projective space  whose topology 
subbasis is its affine covering.
It is
a finite space 
encoding the ``combinatorics" of an $N$-covering, 
 and is non-Hausdorff unless $N=1$. The lattice of all
open subsets of ${{\mathbb P}^{N-1}({\mathbb Z/2})}$ turns out
to be isomorphic to a certain lattice of antichains in the set
of all subsets of an $N$-element set. Combining this with the
Chinese Remainder Theorem for distributive lattices of ideals
in an arbitrary ring, we prove that  distributive
lattices generated by $N$ ideals are equivalent to flabby sheaves 
 over ${{\mathbb P}^{N-1}({\mathbb Z/2})}$.

In particular, consider a compact
Hausdorff space $X$ with a
covering by $N$ closed subsets $X_1,\ldots,X_N$.
Then we have the soft sheaf of continuous functions
with $N$ distinguished $C^*$-algebras $C(X_1),\ldots,C(X_N)$.
However, the soft sheaf of complex-valued continuous functions
on $X$ is not a sheaf of $C^*$-algebras. Therefore, there seems
to be no evident way to use soft sheaves in the noncommutative setting. 
To overcome this difficulty, we declare the closed sets open and consider
$X$ with the new topology generated by these open sets. This leads us
to flabby sheaves over ${{\mathbb P}^{N-1}({\mathbb Z/2})}$. 

This way
we obtain a covering version of Gelfand's theorem: there is 
an equivalence between the category of compact Hausdorff spaces with
$N$-coverings by closed subsets and 
the category of flabby sheaves of unital
commutative $C^*$-algebras over ${{\mathbb P}^{N-1}({\mathbb Z/2})}$. 
In the noncommutative setting, this sheaf-theoretic reformulation 
of coverings allows us to view piecewise trivial comodule algebras 
as introduced
in this paper as what is called ``locally trivial quantum principal
bundles" in \cite{p-mj94}.

Throughout, we work over a field $k$ and all considered
algebras, coalgebras etc.~are over $k$. An unadorned 
$ \otimes $ denotes the tensor product of $k$-vector
spaces. For coproduct and coactions we adopt the Sweedler notation
with the summation sign suppressed: 
$\Delta(h)=h_{(1)} \otimes h_{(2)}\in H\otimes H$,
$\Delta_P(p)=p_{(0)} \otimes p_{(1)}\in P\otimes H$.

\section{Background}\label{prelimi}

\subsection{Fibre products}

We recall here elementary facts concerning pullback diagrams
that will be used in what follows. To focus attention, we
consider the category of vector spaces, which suffices for
our applications.
Let $\pi_1:V_1\rightarrow V_{12}$ 
and $\pi_2:V_2\rightarrow V_{12}$ be linear maps of 
vector spaces. The pullback (fibre product)
$V_1\times_{\pi_1,\pi_2}V_2$
of $\pi_1$ and $\pi_2$ is defined
by a universal property, and turns out to be isomorphic to  
\[  
 \ker\left(\pi_1-\pi_2 :V_1\times V_2\longrightarrow V_{12}\right) 
=\{(p,q) \in V_1 \times V_2\,|\,\pi_1(p)=\pi_2(q)\}.
\] 

As a consequence of this description, we obtain:
\[\label{compa}
		  (\underset{\pi_1,\pi_2}{V_1\times V_2})^{\otimes 2} 
		  =\ker((\pi_1-\pi_2)\otimes\id)\cap\ker(\id\otimes(\pi_1-\pi_2)).
\]
This is because for any linear map 
$ f : V \rightarrow W $ of vector spaces one has 
$ \mathrm{ker}\,f \otimes \mathrm{ker}\, f=
(\mathrm{ker}\, f \otimes V) \cap (V \otimes
\mathrm{ker}\, f) $.

Next, let us consider
the following commutative diagram of linear maps
\begin{equation}\label{hunz}
		  \xymatrix@=8mm{ & & V 
\ar[d]_\eta \ar[ddll]_{\phi_1} \ar[ddrr]^{\phi_2} & &\\ 
& & \overset{\phantom{\pi}}{V_1}\!\underset{\pi_1,\pi_2}{\times} \!V_2 
\ar[dll]^{p_1} \ar[drr]_{p_2} & &\\
V_1 \ar[drr]_{\pi_1} & & & & V_2 \ar[dll]^{\pi_2}\ ,\\
& & V_{12} & &}
\end{equation} 
and show:
\begin{lemma}\label{fib}
Assume that the $\phi_i$'s and $\pi_i$'s in
the  diagram (\ref{hunz}) are surjective.
Then $\eta$ is surjective if and only if 
$\ker\,(\pi_i\!\circ\!\phi_i)=
\mathrm{ker}\,\phi_1+\mathrm{ker}\,\phi_2$.
\end{lemma}
 \begin{proof}
Assume first that $\eta$ is surjective, and 
$v \in \mathrm{ker}\, \pi_i \circ \phi_i$. Then
both $(\phi_1(v),0)$ and 
$(0,\phi_2(v))$ belong to 
$V_1\times_{\pi_1,\pi_2} V_2$, and 
there exist $v_1$ and $v_2$ such that $\eta(v_1)=(\phi_1(v),0)$
and $\eta(v_2)=(0,\phi_2(v))$. Clearly, $v-(v_1+v_2)\in\ker\eta$.
Therefore, as $v_1\in\ker\phi_2$, $v_2\in\ker\phi_1$, and
\[\label{kereta}
\mathrm{ker}\, \eta=
\mathrm{ker}\,\phi_1 \cap
\mathrm{ker}\,\phi_2\,,
\]
we conclude that $v\in \ker\phi_1+\ker\phi_2$.
 
Conversely, if
$\mbox{$(\phi_1(v_1),\phi_2(v_2))$}$ is in the fibre product,
then
$v_1-v_2\in\ker\phi_1+\ker\phi_2 $, so that 
$v_1-v_2=k_1+k_2$ for some $k_1\in\ker\phi_1$
and $k_2\in\ker\phi_2$. 
Hence, for $v:=v_1-k_1=v_2+k_2$, we have 
$(\phi_1(v_1),\phi_2(v_2))=(\phi_1(v),\phi_2(v))=
\eta(v)$. 
\end{proof}

\subsection{Distributive lattices}

We need
a method  yielding a  presentation
of elements in finitely generated distributive
lattices. 
We denote by $\Achain_N$ the set of antichains in $2^{\{1,\ldots,N\}}$.
For any $l\subset 2^{\{1,\ldots,N\}}$ we define
\begin{gather}
\min l\equiv\{u\in l\;|\;\forall v\subsetneq u: v\notin l\},\\
\Upper(l)\equiv \{u\subset\{1,\ldots,N\}\;|\;\exists v\in l:v\subseteq u\}.
\end{gather}
It is easy to see that $\min l\in\Achain_N$.
If $ (\Lambda,\vee,\wedge)$ is a
lattice generated by $\lambda_1,\ldots,\lambda_N$, then we 
define  map $L^\Lambda:\Lambda\rightarrow \Achain_N$, 
\begin{equation}
L^\Lambda(\lambda)=\min\{\{i_1,\ldots,i_k\}\subset \{1,\ldots,N\}\;|\;\lambda_{i_1}\wedge\ldots\wedge\lambda_{i_k}\leq \lambda\} 
\end{equation}
Conversely, we define map $R^\Lambda:2^{2^{\{1,\ldots,N\}}}\rightarrow \Lambda$,  
\begin{equation}\label{stand}
R^\Lambda(l)=\bigvee_{(i_1,\ldots,i_k)\in l}
(\lambda_{i_1} \wedge\ldots\wedge \lambda_{i_k}).
\end{equation}

By the definition of $L^\Lambda$ it follows that for any $\lambda\in\Lambda$, $R^\Lambda(L^\Lambda(\lambda))\le\lambda$. 
On the other hand, as $\Lambda$ is a distributive lattice generated by the
$\lambda_i$'s, there exits $l\subset 2^{\{1,\ldots,N\}}$ such that $\lambda=R^\Lambda(l)$.
Note that for any $\{i_1,\ldots,i_k\}\in l$, $\lambda_{i_1}\wedge\ldots\wedge\lambda_{i_k}\le \lambda=R^\Lambda(l)$, and hence there exists 
$\{j_1,\ldots,j_n\}\in L^\Lambda(\lambda)$ such that 
$\{j_1,\ldots,j_n\}\subseteq \{i_1,\ldots,i_k\}$.
It follows that
$\lambda_{i_1}\wedge\ldots\wedge\lambda_{i_k}\le R^\Lambda(L^\Lambda(\lambda))$ and therefore $\lambda\le R^\Lambda(L^\Lambda(\lambda))$.
Thus, we have proven that for any finitely generated distributive lattice $\Lambda$
\begin{equation}\label{latticeid}
R^\Lambda\circ L^\Lambda=\text{id}_\Lambda.
\end{equation}

Let us define two binary operations $\wedge,\vee:\Achain_N\times\Achain_N\rightarrow \Achain_N$:
\begin{equation}
\label{meetjoin}
l_1\wedge l_2=\min\{u_1\cup u_2\;|\;u_1\in l_1,\; u_2\in l_2\},\ \ \ 
l_1\vee l_2=\min\{l_1\cup l_2\}.
\end{equation}
It is immediate by the distributivity of $\Lambda$ that for all $l_1,l_2\in\Achain_N$,
\begin{equation}\label{mordot}
R^\Lambda(l_1\wedge l_2)=R^\Lambda(l_1)\wedge R^\Lambda(l_2),\ \ \ 
R^\Lambda(l_1\vee l_2)=R^\Lambda(l_1)\vee R^\Lambda(l_2)
\end{equation}

\subsection{Principal extensions} 

Let $(H,\Delta,\varepsilon,S)$ be a Hopf algebra 
with bijective antipode. 
A right $H$-comodule algebra $P$ is a unital associative 
algebra  equipped with an $H$-coaction 
$ \Delta_P : P \rightarrow P \otimes H$ that is an algebra
map.  For a comodule algebra $P$, we call 
\[
	P^\co:=\{p \in P\,|\,\Delta_P(p)=p \otimes 1\}
\] 
the subalgebra of
coaction-invariant elements in $P$.
The assumed existence of the  inverse of the antipode allows us to define 
a left coaction  
${}_P \Delta : P \rightarrow H \otimes P$ 
 by the formula $p \mapsto S^{-1}(p_{(1)}) \otimes p_{(0)}$. This makes
$P$ a left $H$-comodule and a left $H^{op}$-comodule algebra.
\begin{dfn}
Let $P$ be a right comodule algebra over a Hopf algebra 
$H$ with bijective antipode, and let $B:=P^\co$ be the coaction-invariant
subalgebra. The comodule algebra $P$
is called {\em principal} if the following conditions are satisfied:
\begin{enumerate}
\item
the coaction of $H$ is Galois, that is, the map
\vspace*{-1.5mm}$$
\can : P \underset{B}{\otimes} P \longrightarrow P \otimes H,\quad
	p \otimes q \longmapsto pq_{(0)} \otimes q_{(1)}\,,
\vspace*{-1.5mm}$$
(called the canonical map) is bijective,  
\vspace*{3mm}\item
the comodule algebra $P$ is right
$H$-equivariantly projective as a left  $B$-module, i.e.,
 there exists
 a right $H$-colinear and left $B$-linear 
splitting of the multiplication map 
$B \otimes P \rightarrow P$.
\end{enumerate} 
\end{dfn} 
\noindent This splitting can always 
be chosen to be unital
\cite{bh04,bhms+}. Also, in this setting,  
one can show that the $H$-equivariant projectivity of $P$ over $B$
is equivalent to the faithful flatness of $P$ as a $B$-module \cite{ss05,ss}.
If $P$ is a principal comodule algebra, then the extension of algebras
$B\!\subseteq\! P$ is a special case of principal extensions defined in
\cite{bh04}. 

A cleft Hopf-Galois extension (e.g.,  a smash product 
$B \# H$ of an $H$-module algebra $B$ by $H$) is
always principal. Indeed, by \cite[p.42]{bm89}, cleft Hopf-Galois 
extensions
always enjoy the normal basis property, and the latter can be viewed as
equivariant freeness, a special case of equivariant projectivity.
For more details and an introduction to Hopf-Galois theory,
 see, e.g., \cite{m-s93,s-hj94}.

\subsection{Strong connections} 

The inverse of the canonical map defines a monomorphism
$H \rightarrow P \otimes_B P$, 
$h \mapsto \can^{-1}(1 \otimes h)$ called the translation map. 
It turns out that 
lifts of this map to $P \otimes P$ that are both right
and left $H$-colinear yield
an equivalent 
approach to principality \cite{bh04} :
\begin{dfn}\label{strongdef}
Let $H$ be a Hopf algebra with bijective antipode.
Then a strong connection (cf.~\cite{h-pm96,dgh01}) on 
a right $H$-comodule algebra $P$ is a unital linear map
$\ell : H \rightarrow P \otimes P$
satisfying 
$$
(\mathrm{id}_P \otimes \Delta_P) \circ 
\ell = (\ell \otimes \mathrm{id}_H) \circ \Delta,\quad 
({}_P \Delta \otimes \mathrm{id}_P) \circ 
\ell = (\mathrm{id}_H \otimes \ell) \circ
\Delta\,,\quad
\tcan \circ \ell=1 \otimes \id.
$$
Here 
$\tcan : P \otimes P \rightarrow P \otimes H$
is the  lift of $\can$ to $P \otimes P$.
\end{dfn} 
The last property of the strong connection (splitting of $\tcan$)
 gives rise to the 
commutative diagram
\[
\xymatrix{H \ar[r]^-{\ell} \ar[d]_{1\otimes\id} & 
P \otimes P \ar[dl]_{\tcan} \ar[d]^{\text{canonical surjection}}\\
P \otimes H & \;\;P \otimes_B P \,.\ar[l]^{\,\can}}
\]
Using the Sweedler-type notation
$h \mapsto \ell(h)^{\langle 1 \rangle} \otimes 
 \ell(h)^{\langle 2 \rangle}$ (summation suppressed), we can
write the bicolinearity and splitting property of a strong connection 
as follows:
\begin{align}
&\ell(h)^{\langle 1 \rangle}\otimes{\ell(h)^{\langle 2 \rangle}}_{(0)}
\otimes {\ell(h)^{\langle 2 \rangle}}_{(1)}
=\ell(h_{(1)})^{\langle 1 \rangle}\otimes \ell(h_{(1)})^{\langle 2 \rangle}
\otimes h_{(2)}\,,\\
&{\ell(h)^{\langle 1 \rangle}}_{(0)}\otimes {\ell(h)^{\langle 1 \rangle}}_{(1)}
\otimes \ell(h)^{\langle 2 \rangle}
=\ell(h_{(2)})^{\langle 1 \rangle}\otimes S(h_{(1)})\otimes 
\ell(h_{(2)})^{\langle 2 \rangle}\,,\\
&\ell(h)^{\langle 1 \rangle}{\ell(h)^{\langle 2 \rangle}}_{(0)}\otimes
{\ell(h)^{\langle 2 \rangle}}_{(1)}=1\otimes h\,.
\end{align}
Applying $\id\otimes\varepsilon$ to the last equation yields a very useful
formula:
\[\label{elco}
\ell(h)^{\langle 1 \rangle}{\ell(h)^{\langle 2 \rangle}}=\varepsilon(h).
\]

One can prove that an $H$-comodule algebra $P$  is principal
if and only if it admits a strong connection \cite{bh04,bhms+,bb05}.
 Given
a strong connection $\ell$, one can show that the formula 
\[\label{ei0}
P \otimes H \longrightarrow P \underset{B}{\otimes} P,\quad
	p \otimes h \longmapsto p\ell(h)^{\langle 1 \rangle} \otimes 
	\ell(h)^{\langle 2 \rangle},
\] 
defines the inverse of the canonical map $\can$, so that the
coaction of $H$ is Galois. Next, one can also show that
\begin{equation}\label{ei1}
	\s : P\ni p \longmapsto p_{(0)}\ell(p_{(1)})^{\langle 1 \rangle} 
	\otimes \ell(p_{(1)})^{\langle 2 \rangle} \in 
	B \otimes P. 
\end{equation}
is a  splitting whose existence proves the equivariant projectivity.
Much as above, one argues that the formula
\[\label{ei2}
\s' : P\ni p \longmapsto \ell(S^{-1}(p_{(1)}))^{\langle 1 \rangle} 
	\otimes \ell(S^{-1}(p_{(1)}))^{\langle 2 \rangle}p_{(0)} \in 
	P \otimes B 
\]
provides a left $H$-colinear and right $B$-linear splitting of the
multiplication map\linebreak $P\otimes B\rightarrow P$.

\subsection{Actions of compact quantum groups}

The functions continuous along the base and polynomial
along the fibre on a principal  bundle with
compact structure group have an analogue in the
noncommutative setting: Let $\bar H$ be the 
$C^*$-algebra of a compact quantum group in the sense
of Woronowicz \cite{w-sl87,w-sl98} and 
$H$ its dense Hopf $*$-subalgebra
spanned by the matrix coefficients of the irreducible
unitary corepresentations. Let $\bar P$ be a unital $C^*$-algebra 
and let $\delta:\bar P\rightarrow \bar P\otimes_{\mathrm{min}} \bar H$
be a $C^*$-algebraic right
coaction of $\bar H$ on $\bar P$.  
(See \cite[Definition 0.2]{bs93} for a general definition and 
  \cite[Definition 1]{b-fp95} for the special case of compact quantum groups.) 
Then the subalgebra $P\subset\bar P$ of elements for
which the coaction lands in $\bar{P} \otimes H$ 
(algebraic tensor product),
\[
P:=\{p\in \bar P\,|\,\delta(p)\in \bar P\otimes H\},
\]
is an $H$-comodule algebra. It follows from results of 
\cite{b-fp95} and \cite{p-p95} that $P$ is dense in $\bar P$.
We call $P$ the comodule algebra 
associated to the
$C^*$-algebra $\bar P$. 
 It is straightforward to verify that the
operation $\bar P \mapsto P$ commutes with 
taking fibre products.
Note also that $\bar P^{co\bar H}=P^{co H}$.

\section{Piecewise principality}\label{mare}

To show the piecewise nature of principality, we begin by proving
lemmas concerning quotients and fibre products of principal comodule
algebras.

\begin{lemma}\label{fer}
Let $ \pi : P \rightarrow Q$ be a surjection of right $H$-comodule 
algebras (bijective antipode assumed). If $P$ is principal, then:
\vspace*{-2mm}\begin{enumerate}
\item The induced map 
$\pi^\co : P^\co \rightarrow Q^\co$ is surjective.
\item There exists a unital $H$-colinear splitting of $\pi$.
\end{enumerate}  
\end{lemma}
\begin{pf}
It follows from the colinearity of $\pi$ that 
$ \pi (P^\co) \subseteq Q^\co$. To prove the
converse inclusion,  we take advantage of the left $P^\co$-linear
retraction of the inclusion $P^\co\subseteq P$ that was used to
prove \cite[Theorem~2.5(3)]{bh04}:
\[
\sigma:P\longrightarrow P^\co,\quad 
\sigma(p):=p_{(0)}\ell(p_{(1)})^{\langle 1 \rangle}
\varphi(\ell(p_{(1)})^{\langle 2 \rangle})\,.
\]
Here $\ell$ is a strong connection on $P$ and $\varphi$ is any unital
linear functional on $P$. If $\pi(p)\in Q^\co$, then 
$\sigma(p)$ is a desired element of $P^\co$ that is mapped by $\pi$  
to $\pi(p)$. Indeed, since $\pi(p_{(0)})\otimes p_{(1)}=\pi(p)\otimes 1$,
using the unitality of $\ell$, $\pi$ and $\varphi$, we compute
\[
\pi(\sigma(p))=\pi(p_{(0)})\pi(\ell(p_{(1)})^{\langle 1 \rangle})
\varphi(\ell(p_{(1)})^{\langle 2 \rangle})=\pi(p).
\]

Concerning the second assertion, one can readily verify that the formula 
\[
\varsigma(q):= 
\alpha^\co (q_{(0)}\pi (\ell(q_{(1)})^{\langle 1 \rangle})) 
\ell(q_{(1)})^{\langle 2 \rangle}
\]
 defines a unital colinear splitting of $\pi$.
Here $\alpha^\co$ is any $k$-linear unital splitting of $\pi^\co$,
e.g., $\alpha^\co=\varsigma^\co$, and $\ell$ is again a strong connection on $P$.
\end{pf}\\
\begin{lemma}\label{strocopuba}
Let $P$ be a fibre product in the category of right $H$-comodule algebras:
$$
\mbox{$\xymatrix@=5mm{& P \ar[ld] 
\ar[rd] & \\
P_1 \ar[dr]_{\pi^1_2} & & P_2 \ar[dl]^{\pi^2_1}\\
& P_{12}\,. &}$}
$$ 
Then, if $P_1$ and $P_2$ are principal and $ \pi^1_2$ and $\pi^2_1$
are surjective, $P$ is a principal comodule algebra.
\end{lemma} 
\begin{pf}
Given strong connections $\ell_1$ and $\ell_2$ on $P_1$ and $P_2$,
respectively, we want to construct a strong connection on $P$.
A first approximation for such a strong connection is as follows:
\[
\lambda:H\longrightarrow P\otimes P,\quad
\lambda(h):=(\ell_1(h)^{\langle 1 \rangle},f^1_2(\ell_1(h)^{\langle 1 \rangle}))
\otimes (\ell_1(h)^{\langle 2 \rangle},f^1_2(\ell_1(h)^{\langle 2 \rangle})).
\]
Here $f^1_2:=\sigma_2\circ\pi^1_2$, and $\sigma_2$ is a unital colinear
splitting of $\pi^2_1$, which exists by Lemma~\ref{fer}(2).
The map $\lambda$ is unital and bicolinear, but it does not split the
lifted canonical map:
$$
(1,1)\otimes h-\tcan(\lambda(h))=(0,1)\otimes h-
(0,f^1_2(\ell_1(h_{(1)})^{\langle 1 \rangle})
f^1_2(\ell_1(h_{(1)})^{\langle 2 \rangle}))\otimes h_{(2)}\in P_2\otimes H
\,.
$$
Now, let $\tcan_2$ be the lifted canonical map on $P_2\otimes P_2$. 
Applying the splitting of $\tcan_2$ given by $\ell_2$ (cf.~(\ref{ei0}))
to the right hand side of the above equation, gives a correction term for 
$\lambda$:
\begin{align}
T(h):=&\;\ell_2(h)-f^1_2(\ell_1(h_{(1)})^{\langle 1 \rangle})
f^1_2(\ell_1(h_{(1)})^{\langle 2 \rangle})\ell_2(h_{(2)})^{\langle 1 \rangle}
\otimes \ell_2(h_{(2)})^{\langle 2 \rangle}\\
= &\;\left(\varepsilon(h_{(1)})-f^1_2(\ell_1(h_{(1)})^{\langle 1 \rangle})
f^1_2(\ell_1(h_{(1)})^{\langle 2 \rangle})\right)
\ell_2(h_{(2)})^{\langle 1 \rangle}
\otimes \ell_2(h_{(2)})^{\langle 2 \rangle}.\nonumber
\end{align}
This defines a bicolinear map into $P_2\otimes P_2$ which annihilates 1.
Considering $\lambda$ as a map into $(P_1\oplus P_2)^{\otimes 2}$,
we can add these two maps. 
The map $\lambda+ T$ is still unital and bicolinear and splits the lifted
canonical map on $(P_1\oplus P_2)^{\otimes 2}$. Remembering 
(\ref{compa}) and (\ref{elco}), 
it is clear from the formula for $T$ that
to make it taking values in $P\otimes P$ we only need to add the term
\[\label{t'}
T'(h):=\left(\varepsilon(h_{(1)})-f^1_2(\ell_1(h_{(1)})^{\langle 1 \rangle})
f^1_2(\ell_1(h_{(1)})^{\langle 2 \rangle})\right)
\ell_2(h_{(2)})^{\langle 1 \rangle}
\otimes f^2_1(\ell_2(h_{(2)})^{\langle 2 \rangle}).
\]
Much as above, here $f^2_1:=\sigma_1\circ\pi^2_1$ and $\sigma_1$ is a unital 
colinear
splitting of $\pi^1_2$, which exists by Lemma~\ref{fer}(2). The formula (\ref{t'})
defines a bicolinear map into $P_2\otimes P_1$ which annihilates 1.
Since the lifted canonical map on $(P_1\oplus P_2)^{\otimes 2}$ 
vanishes on $P_2\otimes P_1$,
the sum $\lambda+T+T'$ splits the lifted canonical map, takes values in 
$P\otimes P$ and is unital and bicolinear. Thus it is
 as desired a strong connection
on $P$.
\end{pf}\\

Let us now consider a family $\pi_i : P\rightarrow P_i$, $i\in \{1,\ldots,N\}$,
of surjections of right $H$-comodule algebras
 with $\bigcap_{i=1}^N \ker\pi_i=0$. Denote
$J_i:=\ker\pi_i$.
By Lemma~\ref{fib} and formula~\ref{kereta}, for  any $k=1,\ldots,N-1$
there is a fibre-product diagram of right $H$-comodule algebras
\[
\mbox{$\xymatrix@=5mm{
& P/(J_1 \cap \ldots \cap J_{k+1}) \ar[ld]+<3mm,3mm> 
\ar[rd] & \\
\save-<-3mm,0mm>*{P/(J_1 \cap \ldots \cap J_k)} 
\restore \mbox{\huge \  } 
\save <3mm,-14mm>*{} \ar[dr]_(.15){\pi^1_2} \restore
& & P/J_{k+1} \ar[dl]^(.35){\pi^2_1}\\
& P/((J_1\cap\ldots\cap J_k)+J_{k+1}) &.}$}
\]
Assume that all the $P_i\cong P/J_i$ are principal.
Then
Lemma~\ref{strocopuba} implies by an obvious induction 
that  
$P/(J_1\cap\ldots \cap J_k)$ is principal for all
$k=1,\ldots,N$. In particular ($k=N$), $P$ is principal.  
On the other hand, if $P$ is principal then all the $P_i$'s
are principal:
If $\ell: H\rightarrow P\otimes P$ is a strong connection on $P$,
then
\[
(\pi_i\otimes\pi_i)\circ\ell: H\longrightarrow P_i\otimes P_i 
\] 
is a strong connection on $P_i$.
Thus we have proved the following:
\begin{thm}\label{main}
Let $\pi_i : P\rightarrow P_i$, $i\in \{1,\ldots,N\}$,
 be surjections of right $H$-comodule algebras
 such that $\bigcap_{i=1}^N \ker\pi_i=0$.
Then $P$ is principal if and only if all the 
$P_i$'s are principal.
\end{thm}

Our next step is a statement about a relation between the ideals of
a principal comodule algebra $P$ that are also subcomodules, and ideals
in the subalgebra $B$ of coaction-invariant elements. Both sets are
 obviously
lattices with respect to the operations $+$ and $\cap$.

\begin{prop}\label{diens}
Let $P$ be a principal right  $H$-comodule algebra  and 
$B:=P^\co$ the coaction-invariant subalgebra. Denote by
$\Xi_B$ the lattice of all ideals in $B$
and by $\Xi_P$ the lattice of all ideals  
in $P$ which are simultaneously subcomodules.
Then the map
$$
	\mathcal{L}:\Xi_P \longrightarrow \Xi_B,\quad 
	\mathcal{L}(J) := J \cap B
$$
is a monomorphism of lattices. 
\end{prop}
\begin{pf}
The only non-trivial step in proving that $\mathcal{L}$ 
is a homomorphism of lattices
is establishing the inclusion
$(B \cap J) + (B \cap J') \supseteq B \cap (J+J')$. To this end,
we  proceed along the lines of the proof of Lemma~\ref{fer}(1).
Since $J$ is a comodule and an ideal, from the formula (\ref{ei1}) we obtain:
\[\label{stropro}
p\in J\implies s(p)=p_{(0)}\ell(p_{(1)})^{\langle 1 \rangle}
\otimes \ell(p_{(1)})^{\langle 2 \rangle}\in (J\cap B)\otimes P.
\]
Now, let $p \in J$, $q \in J'$, $p+q \in B$. Then 
(\ref{stropro}) implies that 
\[
	(p+q) \otimes 1 = s(p)+s(q)
\in 
	(B \cap J) \otimes P +
	(B \cap J') \otimes P.
\] 
Applying  
any unital linear functional $P \rightarrow k$
to the second tensor component implies
$p+q \in (B \cap J)+(B \cap J')$. 
Finally, since $s$ is a splitting of the multiplication map,
it follows from 
(\ref{stropro}) that 
$J=(J\cap B)P$.
 This, in turn, proves
the injectivity of $\mathcal L$.
\end{pf}\\

\begin{rem}\rm
Much as the formula (\ref{ei1}) 
implies (\ref{stropro}), the formula (\ref{ei2})
implies
\[\label{stropro2}
p\in J\implies s'(p)=\ell(S^{-1}(p_{(1)}))^{\langle 1 \rangle}
\otimes \ell(S^{-1}(p_{(1)}))^{\langle 2 \rangle}p_{(0)}\in P\otimes
(J\cap B).
\]
 Therefore, since $s'$ is a splitting of the right multiplication map,
$J=P(J\cap B)$. Combining this with the above discussed left-sided 
version yields
\[
\label{linv}
P(J\cap B)=J=(J\cap B)P.
\]
\end{rem}\smallskip

\begin{rem}\rm
Note that the homomorphism $\mathcal L$ is not
surjective in general. A counterexample is given 
by a smash product (trivial principal comodule algebra) of 
the Laurent polynomials $B=k[u,u^{-1}]$ with the Hopf
algebra $H=k[v,v^{-1}]$ of Laurent polynomials 
($ \Delta(v)=v \otimes v$). The action is defined by 
 $v \triangleright u=qu$, 
$q \in k \setminus\! \{0,1\}$. Viewing $u$ and $v$ as generators
of $P$, it is clear that  
	$vu=quv$.
It is straightforward to verify that,  if $I$ 
is the two-sided ideal in $B$ generated by $u-1$, then the
right ideal $IP$ is not a 
two-sided ideal of $P$. Hence the map $\mathcal L$ cannot be surjective by (\ref{linv}).
\end{rem}

In \cite{cm00}, families $ \pi_i : P \rightarrow P_i$ of 
algebra homomorphisms as in Theorem~\ref{main}
were called coverings. However, it was explained therein that 
such coverings are well-behaved
when the kernels $ \mathrm{ker}\, \pi_i$ generate a
distributive lattice of ideals (with respect to $+$ and
$\cap$ as lattice operations). Hence we adopt in the
present paper the following terminology:
\begin{dfn}\label{covdef}
A finite family 
$\pi_i:P\rightarrow P_i$, $i=1,\ldots,N$, of surjective 
algebra homomorphisms
is called a weak covering if 
$\cap_{i=1,\ldots,N}\ker\pi_i=\{0\}$.
A weak covering is 
called a {\em covering} if the lattice of ideals generated by
the $\ker\pi_i$'s is distributive.
\end{dfn}
\noindent
The above definition can  obviously be extended to the
case when the $\pi_i$'s 
are algebra and $H$-comodule morphisms. Then the $\ker\pi_i$'s are 
ideals and $H$-subcomodules.

The next claim is concerned with
the distributivity condition from
Definition~\ref{covdef} for coverings of principal
comodule algebras. It
 follows  from Proposition~\ref{diens} and Lemma~\ref{fer}(1).

\begin{cor}\label{auchmain}
Let  
$ \pi_i : P \rightarrow P_i$, $i=1,\ldots,N$, be surjective
 homomorphisms of right $H$-comodule algebras. Assume that
$P$ is principal. Then
$\{\pi_i : P \rightarrow P_i\}_i$ is a covering of $P$  if and
only if $\{\pi_i^\co : P^\co \rightarrow P_i^\co\}$ is a covering of 
$P^\co$.
\end{cor}
\noindent
The above corollary  is particularly helpful when $P^\co$ is a
$C^*$-algebra because lattices of closed ideals in a $C^*$-algebra
are always distributive (due to the property $I\cap J=IJ$).

We are now ready to propose a noncommutative-geometric replacement
of the concept of local triviality of principal bundles.
Since these are the closed 
rather than open subsets of a compact Hausdorff
space that admit a natural translation into the language of 
$C^*$-algebras, we use finite closed rather than open coverings
to trivialize  bundles. As is explained in \cite[Example~1.24]{bhms},
there is a difference between these two approaches.
We reserve the term ``locally trivial" for bundles trivializable
over an open cover, and call bundles trivializable over a finite
closed cover ``piecewise trivial". It is the latter (slightly more
general) property that we generalize to the noncommutative setting.
\begin{dfn}\label{pipri}
An $H$-comodule algebra $P$ 
is called \emph{piecewise principal (trivial)} if there exist
comodule algebra surjections 
$ \pi_i : P \rightarrow P_i$, $i=1,\ldots,N$, 
such that:
\begin{enumerate}
\item The restrictions  
$\pi_i|_{P^\co} : P^\co \rightarrow P_i^\co$ form a covering. 
\item The $P_i$'s are principal. (The 
$P_i$'s are isomorphic as 
$H$-comodule algebras to a smashed product 
$P_i^\co \#_i H$.) 
\end{enumerate} 
\end{dfn}  

While not every compact principal bundle is   piecewise
(or locally)
trivial (\cite[Example~1.22]{bhms}), 
every piecewise principal compact $G$-space 
(i.e., covered by finitely many compact principal $G$-bundles) is clearly
a compact principal $G$-bundle.  The second statement becomes non-trivial
when we replace compact $G$-spaces by comodule algebras. However, it is
an immediate consequence of Theorem~\ref{main} and 
Corollary~\ref{auchmain}:
\begin{cor}
Let $H$ be a Hopf algebra with bijective antipode and 
$P$ be an $H$-comodule algebra that is piecewise
principal  with respect to $ \{\pi_i : P \rightarrow P_i\}_i$. 
Then $P$ is principal and $ \{\pi_i : P \rightarrow P_i\}_i$
is a covering of $P$. 
\end{cor}

Finally, let us consider  the relationship
between piecewise triviality  and a similar concept referred
to as ``local triviality" in \cite{p-mj94}. Therein, sheaves
$\P$ of comodule algebras were viewed as quantum
analogues of principal  bundles. They were called 
locally trivial provided that the space $X$ on
which $\P$ is defined admits an open covering 
$\{U_i\}_i$ such that all $\P(U_i)$'s are smash products.  
If we assume such a sheaf to be flabby (that is,  for all
open subsets of $V,U$, $V\subseteq U$, 
the restriction maps 
$\pi_{U,V} : \P(U) \rightarrow \P(V)$
 are surjective),
then we can use 
Theorem~\ref{main}  to deduce the principality 
of all $\P(U)$'s:
\begin{cor}\label{niema}
Let $H$ be a Hopf algebra with bijective antipode and
$\P$ be a flabby sheaf of $H$-comodule algebras
over a topological
space $X$. If $\{U_i\}_i$ is a finite open
covering such that all 
$\P(U_i)$'s are principal, then 
$\P(U)$ is principal for any open subset $U\subseteq X$.
\end{cor}

\section{Coverings and flabby sheaves}\label{corollaries}

In this section, we focus entirely on a flabby-sheaf interpretation
of distributive lattices of ideals defining coverings of algebras 
(see Definition~\ref{covdef}).
We will explain
 that for a flabby sheaf in Corollary~\ref{niema}, the underlying
topological space plays only a secondary role and can
be replaced by a certain space that is universal for all $N$-element
coverings. This space is  the 2-element field $N-1$-projective space  
\begin{equation}
{{\mathbb P}^{N-1}({\mathbb Z/2})}:=\{0,1\}^N \setminus\! \{(0,\ldots,0)\}
\end{equation} 
whose topology 
subbasis is its affine covering, i.e., 
it is the topology generated by 
the subsets
\begin{equation} \label{AI}
A_i:=\{(z_1,\ldots,z_N) 
\in {{\mathbb P}^{N-1}({\mathbb Z/2})}\,|\,z_i \neq 0\}.
\end{equation} 
 
Consider now  an arbitrary space $X$
with a finite covering $\{U_1,\ldots,U_N\}$.
Define on $X$
the topology generated by the $U_i$'s (considered
as open sets) and pass to the quotient by 
the equivalence relation 
\begin{equation}\label{paul}
x \sim y \Leftrightarrow (\forall i : 
x \in U_i \Leftrightarrow y \in U_i). 
\end{equation} 
Obviously, 
${X/\!\!\sim}$ depends on the specific features of the
covering $\{U_i\}_i$. However,
for a fixed $N$, it can always be 
embedded into $ {{\mathbb P}^{N-1}({\mathbb Z/2})}$:
\begin{prop}\label{embedding}
\newcommand{\mquos}{{X/\!\sim}}
Let $X=U_1 \cup \ldots \cup U_N$ be any set
equipped with the topology
 generated by the $U_i$'s.
Let $p:X\rightarrow\mquos$ be the 
quotient map defined by the equivalence relation (\ref{paul}). Then
\begin{equation*}
\xi:\mquos\longrightarrow {{\mathbb P}^{N-1}({\mathbb Z/2})},\quad 
p(x)\longmapsto (z_1,\ldots, z_N),\quad
z_i=1\Leftrightarrow x\in U_i\,,\;\forall\, i, 
\end{equation*}
is an embedding of topological spaces.
\end{prop}
\begin{pf}
\newcommand{\mquos}{{X/\!\sim}}
It is immediate  that $\xi$ is well
defined and injective. Next, since $p^{-1}(\xi^{-1}(A_i))=U_i$ is
open for each $i$, all $ \xi^{-1}(A_i)$'s are open   
in the quotient topology on $\mquos$.
Now the  
continuity of $\xi$ follows from the fact that 
$A_i$'s form a subbasis of the topology
of $ {{\mathbb P}^{N-1}({\mathbb Z/2})}$. 

The key step is to show  that images of 
open sets in $\mquos$ are open in $\xi(X/\!\sim)$.
First note that by the definition of the relation~(\ref{paul}), 
\[
p^{-1}(p(U_{i_1}\cap\ldots\cap U_{i_n}))=
U_{i_1}\cap\ldots\cap U_{i_n}\,.
\] 
Therefore, as preimages and images preserve 
  unions and  any open set in $X$ is the union of
  intersections of $U_i$'s, $p$ is an open map.
On the other hand, by the surjectivity of $p$, we have
  $p(p^{-1}(V))=V$ for any subset $V\subset\mquos$.  
  Hence  it follows that 
  a set in $\mquos$ is open if and only if it is an 
image under $p$ of an open set in $X$.
Finally,  by the definition of $\xi$,
\[
\xi(p(U_{i_1}\cap\ldots\cap U_{i_n}))
=A_{i_1}\cap\ldots\cap A_{i_n}\cap\text{Im}(\xi),
\]
 and the claim follows from the distributivity of $\cap$ with respect 
to $\cup$. 
 \end{pf}\\

Note that the map $ \xi $ is a homeomorphism 
precisely when the $U_i$'s are in a generic position, 
that is when all intersections 
$U_{i_1} \cap \ldots \cap U_{i_k}\cap (X\setminus U_{j_1})\cap\ldots\cap
(X\setminus U_{j_l})$ such that 
$\{i_1,\ldots,i_k\}\cap\{j_1,\ldots,j_l\}=\emptyset$ are
non-empty.

Thus we have shown that 
if we consider $X$  and 
the $U_i$'s as in Corollary~\ref{niema}, then the composition 
$ \xi \circ p : X \rightarrow {{\mathbb P}^{N-1}({\mathbb Z/2})}$ 
is continuous.
Hence we can produce 
flabby sheaves over $ {{\mathbb P}^{N-1}({\mathbb Z/2})}$ 
by taking direct
images 
of flabby sheaves over $X$. 
They will have the 
same sections globally and on the covering
sets. In this sense, they carry an essential part of
the data encoded in the original sheaf.

\begin{mex}
\rm
Let $X={\mathbb P}^{N-1}({\mathbb C})$.
Denote by $[x_1 : \ldots : x_{N}]$ the class of 
$(x_1,\ldots,x_{N})\in{\mathbb C}^{N}$
in ${\mathbb P}^{N-1}({\mathbb C})$. 
Define a family of  closed sets
\begin{equation*}
X_i=\{[x_1:\ldots:x_N]\in {\mathbb P}^{N-1}({\mathbb C})\;|\;
|x_i|=\max(\{|x_1|,\ldots,|x_{N}|\})\;\},\ \ \ i=1,\ldots,N.
\end{equation*}
The $X_i$'s cover $X$, and moreover, for all 
$\Lambda,\Gamma\subset \{1,\ldots,N\}$ 
such that $ \Lambda \neq \emptyset$,
$\Lambda\cap\Gamma=\emptyset$,
the element
\begin{equation*}
[x_1 : \ldots : x_{N}]\in X, \text{\ where\ }
x_i=\left\{{{\!\!\!1 \text{\ if\ } i\in\Lambda}\atop {0 \text{\ otherwise}}}\right.
\end{equation*} 
belongs to $\bigcap_{i\in \Lambda}X_i\cap\bigcap_{j\in\Gamma}(X\setminus X_j)$.
It follows that the $X_i$'s are in generic position and the map
$\xi:{X/\!\sim}\rightarrow {{\mathbb P}^{N-1}({\mathbb Z/2})}$ 
(Proposition~\ref{embedding}) is a homeomorphism
(if ${X/\!\sim}$ is considered with finite topology as in
Proposition~\ref{embedding}).
Let
\begin{equation*}
V_i=\{[x_1 : \ldots : x_{N}]\in{\mathbb P}^{N-1}({\mathbb C})\;|\;
x_i\neq 0\;\},\ \ \ i=1,\ldots,N
\end{equation*} 
be the standard (open) affine cover of ${\mathbb P}^{N-1}({\mathbb C})$ and $\Psi_i$ be the
homeomorphism
\begin{equation*}
\Psi_i:V_i\rightarrow {\mathbb C}^{N-1},\ [x_1 :\ldots:x_{N}]
\mapsto (x_1/x_i,\ldots,\widehat{x_i/x_i},\ldots,x_{N}/x_i). 
\end{equation*}
Observe that for all $i\in\{1,\ldots,N\}$,
$X_i\subset V_i$ and
\begin{equation*}
\Psi_i(X_i)=\{(y_1,\ldots,y_{N-1})\in{\mathbb C}^{N-1}\;|\;\forall_j\;
|y_j|\leq 1\;\}.
\end{equation*}
In particular, $X_i$'s are indeed closed sets.
\end{mex}

Our next aim is to demonstrate 
that flabby sheaves over $ {{\mathbb P}^{N-1}({\mathbb Z/2})}$
are just a reformulation of the notion of
covering introduced in Definition~\ref{covdef}. It turns out that the
distributivity condition discussed in the previous
section is the key property needed to reconcile the
results from
\cite{bk96,cm02} with those from \cite{p-mj94}. 
 In particular, our results imply the principality of 
Pflaum's noncommutative instanton bundle 
(see the last section for details).


Let us consider the lattice $\Gamma_N$ of open subsets of ${{\mathbb P}^{N-1}({\mathbb Z/2})}$ generated by the $A_i$'s (\ref{AI}).
\begin{lemma}\label{upperin}
If $A_{i_1}\cap\ldots\cap A_{i_k}\subseteq R^{\Gamma_N}(l)$, for some $l\subset 2^{\{1,\ldots,N\}}$ then $\{i_1,\ldots,i_k\}\in\Upper(l)$.
\end{lemma}
\begin{pf}
Suppose that $\{i_1,\ldots,i_k\}\notin\Upper(l)$. Therefore, for all $u\in l$ there exists $j(u)\in u$, such that $j(u)\notin \{i_1,\ldots,i_k\}$.
Then
\begin{multline*}
A_{i_1}\cap\ldots\cap A_{i_k}\cap(\Proj \setminus R^{\Gamma_N}(l))
=A_{i_1}\cap\ldots\cap A_{i_k}\cap\bigcap_{u\in l}\bigcup_{j\in u}
(\Proj\setminus A_j)\\
\supseteq A_{i_1}\cap\ldots\cap A_{i_k}\cap\bigcap_{u\in l}(\Proj\setminus A_{j(u)})\neq\emptyset,
\end{multline*}
as $A_i$'s are in generic position (see the remark after Proposition~\ref{embedding}). Hence
$A_{i_1}\cap\ldots\cap A_{i_k}$ cannot be contained in $R^{\Gamma_N}(l)$.
\end{pf}

By Lemma~\ref{upperin}, it is clear that
\begin{equation}
\{\{i_1,\ldots,i_k\}\subset\{1,\ldots,N\}\;|\;A_{i_1}\cap\ldots\cap A_{i_k}\subseteq R^{\Gamma_N}(l)\}=\Upper(l).
\end{equation}
As $\min\Upper(l)=\min l$, it follows that for any $l\subset 2^{\{1,\ldots,N\}}$,
\begin{equation}
\label{idcone}
L^{\Gamma_N}(R^{\Gamma_N}(l))=\min l.
\end{equation}
Thus, for any $U,U'\subset {\mathbb P}^{N-1}({\mathbb Z/2})$, using (\ref{idcone}), (\ref{latticeid}) and (\ref{mordot}),
\begin{multline}
L^{\Gamma_N}(U\cap U')=L^{\Gamma_N}(R^{\Gamma_N}(L^{\Gamma_N}(U))\cap R^{\Gamma_N}(L^{\Gamma_N}(U')))\\
=L^{\Gamma_N}(R^{\Gamma_N}(L^{\Gamma_N}(U)\wedge L^{\Gamma_N}(U')))
=\min (L^{\Gamma_N}(U)\wedge L^{\Gamma_N}(U'))
=L^{\Gamma_N}(U)\wedge L^{\Gamma_N}(U').
\end{multline}
Similarily one can show that
\begin{equation}
L^{\Gamma_N}(U\cup U')=L^{\Gamma_N}(U)\vee L^{\Gamma_N}(U').
\end{equation}
Hence, we have proven
\begin{lemma}\label{isomlattice}
The map $L^{\Gamma_N}:\Gamma_N\rightarrow \Achain_N$ is an isomorphism of distributive lattices (with 
$(L^{\Gamma_N})^{-1}=\left.R^{\Gamma_N}\right|_{\Achain_N}$) where $\Achain_N$ is a distributive lattice with meet and join operations defined in
(\ref{meetjoin}).
\end{lemma}

\newcommand{\Kpi}{{\Lambda_{(\ker\pi_i)_i}}}
\newcommand{\Keta}{{\Lambda_{(\ker\eta_i)_i}}}

In the following proposition, 
we use the above results for 
the lattice $\Gamma_N$ of open subsets 
of ${{\mathbb P}^{N-1}({\mathbb Z/2})}$  and 
the lattice $\Kpi$ of ideals
 generated by the kernels of surjections
forming an $N$-covering. For the first
lattice, we consider the category of flabby sheaves over 
${{\mathbb P}^{N-1}({\mathbb Z/2})}$, and for the second, 
the category of $N$-coverings of  
algebras. Here 
``$N$-covering'' means a covering by $N$ surjections, and
a morphism between $N$-coverings
$\{\pi_i : P \rightarrow P_i\}_i$ and
$\{\eta_i : Q \rightarrow Q_i\}_i$ 
consists of morphisms 
$ \xi : P \rightarrow Q$ and 
$ \xi_i : P_i \rightarrow Q_i$ such that 
$ \eta_i \circ \xi = \xi_i \circ \pi_i\,,\forall\, i\in\{1,\ldots,N\}$.
\begin{thm}\label{nagut} 
Let $\mathbf{C}_N$ be the category 
of $N$-coverings of algebras, and
$\mathbf{F}_N$ be 
the category of flabby sheaves of algebras over
${{\mathbb P}^{N-1}({\mathbb Z/2})}$. Then 
the following assignments
\begin{gather}
\label{ggb}
		  \mathbf{C}_N \ni
		  \{\pi_i : P \rightarrow P_i\}_i \longmapsto 
		  \{\P : U \mapsto P/R^{\Kpi}(L^{\Gamma_N}(U))\}_U \in
		  \mathbf{F}_N,\\
\label{gga}
		  \mathbf{F}_N \ni
	\P  \longmapsto 
	\{\P({{\mathbb P}^{N-1}({\mathbb Z/2})}) \rightarrow 
	\P(A_i)\}_i \in \mathbf{C}_N
\end{gather} 
yield an equivalence of categories.
\end{thm}
\begin{pf}
Suppose we are given a 
flabby sheaf $\P$ of algebras over ${{\mathbb P}^{N-1}({\mathbb Z/2})}$.
Let $\pi_{V,U}:\P(V)\rightarrow\P(U)$ denote the restriction map for any 
open $U$,$V$, $U\subseteq V$. 
For brevity, we write $\pi_U$ instead of $\pi_{V,U}\,$, if 
$V={{\mathbb P}^{N-1}({\mathbb Z/2})}$. 
By the flabbiness of $\P$, the morphisms $\pi_{A_i}$ (see (\ref{AI}))
 are surjective. 
The property $\bigcap_{i=1}^N \mathrm{ker}\, \pi_{A_i}=\{0\}$
 follows from the sheaf condition. 

It remains to prove the distributivity of
the lattice  generated by the kernels of 
$\pi_{A_i}$'s. 
Lattices of sets are always distributive, so that
 it is enough to show that the assignment 
$U\mapsto\mathrm{ker}\, \pi_{U}$ defines a surjective morphism
from the lattice of open subsets of 
${{\mathbb P}^{N-1}({\mathbb Z/2})}$ onto the lattice of ideals generated
by $\mathrm{ker}\, \pi_{A_i}$'s. Here by a morphism of lattices we
mean a map that
transforms the union and intersection of open subsets to
the intersection and sum of ideals, respectively. 

To show this, let 
$U'$, $U''$ 
be open subsets of ${\mathbb P}^{N-1}({\mathbb Z/2})$.
Since $\P$ is a sheaf, we know that $\P(U' \cup U'')$ is the
 fibre product
of $\P(U')$ and $\P(U'')$. Now it follows from (\ref{kereta}) that
$
\mathrm{ker}\, \pi_{U' \cup U''}=
\mathrm{ker}\, \pi_{U'} \cap \mathrm{ker}\, \pi_{U''}
$, as needed.
Similarly, since the sheaf $\P$ is flabby,
Lemma~\ref{fib} implies that
$\mathrm{ker}\, \pi_{U' \cap U''}=
\mathrm{ker}\, \pi_{U'} +
\mathrm{ker}\, \pi_{U''}$.
Thus we have shown that (\ref{gga}) assigns coverings to flabby sheaves.

Conversely, assume that we are given a covering.
For brevity, let us denote for any $U\in\Gamma_N$, $\widehat{L(U)}:=R^\Kpi(L^{\Gamma_N}(U))$.
Let $U$, $U'$ be open subsets of ${{\mathbb P}^{N-1}({\mathbb Z/2})}$ such that $U\subseteq U'$.
For all $\{i_1,\ldots,i_k\}\in L^{\Gamma_N}(U)$, $A_{i_1}\cap\ldots\cap A_{i_k}\subseteq U\subseteq U'$ which implies that 
there exists a subset $u\subset\{i_1,\ldots,i_k\}$ such that $u\in L^{\Gamma_N}(U')$.
It is then clear that
\begin{multline}
\widehat{L(U')}=\bigcap_{(j_1,\ldots,j_n)\in L^{\Gamma_N}(U')}
(\ker\pi_{j_1}+\ldots+\ker\pi_{j_n})\\ \subseteq
\bigcap_{(i_1,\ldots,i_k)\in L^{\Gamma_N}(U)}
(\ker\pi_{i_1}+\ldots+\ker\pi_{i_k})=\widehat{L(U)},
\end{multline}
and one can define restriction map
\begin{equation}
\pi_{U',U}:\P(U')\rightarrow \P(U),\ p+\widehat{L(U')}\mapsto p+\widehat{L(U)}.
\end{equation} 
Hence $\P$ is a presheaf.

Let $U$ be an open subset of ${{\mathbb P}^{N-1}({\mathbb Z/2})}$ and let $(U_i)_i$ be an open covering of $U$
(i.e. $\bigcup_iU_i=U$).
Suppose that $(p_{U_i})_i$ is a family of elements, where for each $i,j$, $p_{U_i}\in \P(U_i)$, and
$\pi_{U_i, U_i\cap U_j}(p_{U_i})=\pi_{U_j, U_i\cap U_j}(p_{U_j})$. By the distributivity of lattice of ideals generated by $\ker\pi_i$ and 
generalised Chinese remainder theorem, (see 
e.g.~\cite{ss58}, Theorem~18 on p.~280) there exists an element $p_U\in\P(U)$ such that, for all $i$,
$\pi_{U, U_i}(p_U)=p_{U_i}$.
It is easy to see that if $\widehat{L(U\cup U')}=\widehat{L(U)}\cap\widehat{L(U')}$, for any open subsets $U$, $U'$  of ${{\mathbb P}^{N-1}({\mathbb Z/2})}$
then this element is unique. But this follows by Lemma~\ref{isomlattice} and (\ref{mordot}).

Denote assignement in (\ref{ggb}) by $F$ and functor defined by (\ref{gga}) by $G$. The functoriality of $G$ is immediate.
Let $\xi:P\rightarrow Q$, $(\xi_i:P_i\rightarrow Q_i)_i$ be the morphism of $N$-coverings $\{\pi_i : P \rightarrow P_i\}_i$ and
$\{\eta_i : Q \rightarrow Q_i\}_i$. Note that from
$ \eta_i \circ \xi = \xi_i \circ \pi_i\,,\forall\, i\in\{1,\ldots,N\}$ it follows that for all $i$, $\xi(\ker\pi_i)\subseteq \ker\eta_i$.
Accordingly, for any open subset $U$ of ${{\mathbb P}^{N-1}({\mathbb Z/2})}$, 
$\xi(R^\Kpi(L^{\Gamma_N}(U)))\subset R^\Keta(L^{\Gamma_N}(U))$, 
and therefore
for all open $U$, the following family defines a map of corresponding sheaves
\begin{gather}
P/R^\Kpi(L^{\Gamma_N}(U))\longrightarrow Q/R^\Keta(L^{\Gamma_N}(U))\nonumber \\
p+R^\Kpi(L^{\Gamma_N}(U))\mapsto\xi(p)+R^\Keta(L^{\Gamma_N}(U))
\end{gather}

That $G\circ F$ is a functor naturally isomorphic to identity functor in the category of $N$-coverings 
is obvious, as $\widehat{L(A_i)}=\ker\pi_i$ and $P_i\simeq P/\ker\pi_i=\P(A_i)$.

On the other hand suppose that we are given a flabby sheaf $\P$. By the sheaf property, for any open sets $U$, $U'$,
$\P(U\cup U')$ is a fibre product of $\P(U)$ and $\P(U')$ over $\P(U\cap U')$. Then using flabbines, we can apply Lemma~\ref{fib} and formula
\ref{kereta} to conclude that  
\begin{equation}
\ker\pi_{U\cup U'}=\ker\pi_{U}\cap\ker\pi_{U'},\ \ker\pi_{U\cap U'}=\ker\pi_{U}+\ker\pi_{U'}.
\end{equation}
Then by the property that assignement $U\mapsto\widehat{L(U)}$ is a morphism of lattices and again using flabbiness one sees that for all open $U$,
$\P(U)\simeq\P({{\mathbb P}^{N-1}({\mathbb Z/2})})/R^{\Lambda_{(\ker\pi_{A_i})_i}}(L^{\Gamma_N}(U))$, 
But this immedietaly shows that $F\circ G$ is naturally isomorphic to identity functor on $\mathbf{F}_N$, which ends the proof.
\end{pf}\\


Since the intersection of closed ideals in a 
$C^*$-algebra equals their product, lattices of closed
ideals in $C^*$-algebras are always distributive.
Thus we immediately obtain:
\begin{cor}\label{claco}
Compact Hausdorff spaces with a fixed covering by $N$ closed
subsets are equivalent to flabby sheaves of commutative unital
$C^*$-algebras over ${{\mathbb P}^{N-1}({\mathbb Z/2})}$. 
\end{cor}

Here, the zero algebra is allowed
as a unital $C^*$-algebra. This is needed  
if the closed sets are not in generic position. The 
set of unital morphisms from the zero 
$C^*$ algebra to any other is understood to be empty.

\section{Examples}\label{examples}

In this last section we recall from
\cite{dhh,hms06,bhms05,hms06m}
the construction of examples for the above concepts
that illustrate possible areas of applications.

\subsection{A noncommutative join construction} 
If $G$ is a compact group, then the join 
$G * G$ becomes a $G$-principal fibre bundle over the 
unreduced suspension $ \Sigma G $ of $G$, see 
e.g.~\cite{b-ge93}, Proposition~{VII.8.8} or
\cite{bhms}. 
For example, one can obtain the
Hopf fibrations $S^7 \rightarrow S^4$ and 
$S^3 \rightarrow S^2$ in this way using 
$G=SU(2)$ and $G=U(1)$, respectively. 
Recall that $G * G$ is obtained from 
$[0,1] \times G \times G$ by shrinking to a point 
one factor $G$ at $0 \in [0,1]$ and the other factor $G$ 
at $1$.
 
\vspace*{-0mm}\begin{figure}[h]
\[
\includegraphics[width=90mm]{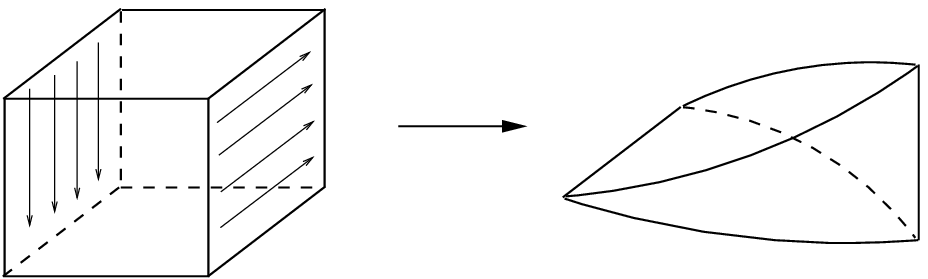}
\nonumber
\]
\end{figure}\vspace*{-0mm}

Alternatively, one can shrink 
$G \times G$ at $0$ to the diagonal. This
picture is generalised in \cite{dhh}.  
Our aim in this
first part of Section~\ref{examples} is to describe a
noncommutative analogue of this construction 
that nicely fits into our
general concepts and will be studied in greater detail
in \cite{dhh}.

To this end, let
$H$ be the Hopf algebra underlying a 
compact quantum group $\bar H$
(see \cite{w-sl87,w-sl98} or Chapter~11 of \cite{ks97} for details). 
We define 
\begin{eqnarray}
&& P_1:=\{f \in C([0,1],\bar H) \otimes H\,|\,
	f(0) \in \Delta(H)\},\nonumber\\ 
&& P_2:=\{f \in C([0,1],\bar H) \otimes H\,|\,
	f(1) \in \mathbb{C} \otimes H\}\nonumber
\end{eqnarray} 
which will play the roles of the
two trivial pieces of the principal extension. Here 
we identify elements of 
$C([0,1],\bar H) \otimes H$ with functions 
$[0,1] \rightarrow \bar H \otimes H$. The $P_i$'s become
$H$-comodule algebras by applying the coproduct 
of $H$ to $H$, 
$ \Delta_{P_i}=\mathrm{id}_{C([0,1],\bar H)} \otimes \Delta$, 
and the subalgebras of $H$-invariants can be identified
with
\begin{eqnarray}
&& B_1:=\{f \in C([0,1],\bar H)\,|\,f(0) \in \mathbb{C}\},
		  \nonumber\\ 
&& B_2:=\{f \in C([0,1],\bar H)\,|\,f(0) \in \mathbb{C}\}.
		  \nonumber 
\end{eqnarray}
Furthermore, 
$P_1 \simeq B_1 \# H$, 
$P_2 \simeq B_2 \otimes H$, where $H$ acts on $B_1$ via
the adjoint action, 
$(a \triangleright f)(t)=a_{(1)}f(t)S(a_{(2)})$,
$a \in H,f \in B_1,t \in [0,1]$,
see \cite{dhh}. Now one can define $P$ as a glueing of
the two pieces along $P_{12}:=\bar H \otimes H$, 
that is, as the pull-back 
$$
	P:=\{(p,q) \in P_1 \oplus P_2\,|\,
	\pi^1_2(p)=\pi_1^2(q)\}
$$
of the $P_i$'s along the
evaluation maps 
$$
	\pi^1_2 : P_1 \rightarrow P_{12},\quad f \mapsto f(1),\quad
	\pi^2_1 : P_2 \rightarrow P_{12},\quad f \mapsto f(0).
$$
Theorem~\ref{main} implies that $P$ is principal. 

\subsection{The Heegaard-type quantum 3-sphere} 
\newcommand{\tplz}{\mathcal{T}}
Based on the idea of a Heegaard splitting of $S^3$ into
two solid tori, a noncommutative deformation of $S^3$
was proposed in \cite{cm02,hms06,bhms05}. 
On the level of $C^*$-algebras, it can be presented as
a fibre product $C(S^3_{pq \theta})$
of two $C^*$-algebraic crossed products
$\tplz \rtimes_{\theta} \mathbb{Z}$ and  
$\tplz \rtimes_{-\theta} \mathbb{Z}$
of the Toeplitz algebra $\tplz$ by $\mathbb{Z}$. 
We denote the isometries generating $\tplz$ in
the two algebras by $z_+,z_-$.  
The $\mathbb{Z}$-actions
are implemented by unitaries 
$u_+,u_-$, respectively, in the following way:
$$		  
u_+ \triangleright_{\theta} z_+=u_+z_+u_+^{-1}:=
e^{2 \pi i \theta} z_+,\quad 
u_- \triangleright_{-\theta} z_-=u_-z_-u_-^{-1}:=
e^{-2 \pi i \theta} z_-.
$$ 
The fibre product is taken over 
$C(S^1) \rtimes_{\theta} \mathbb{Z}$ 
with action
$U_+ \triangleright_{\theta} Z_+:=e^{2 \pi i \theta}Z_+$,
where
$Z_+$ is the generator of $C(S^1)$ and $U_+$ is the
unitary giving the $\mathbb{Z}$-action in this
algebra. 
The corresponding surjections defining the
fibre product are
\begin{eqnarray}
&& \pi^1_2 : \tplz \rtimes_{\theta} \mathbb{Z} \rightarrow 
C(S^1) \rtimes_{\theta} \mathbb{Z},\quad 
z_+ \mapsto Z_+,\quad u_+ \mapsto U_+,\nonumber\\ 
&& \pi^2_1 : \tplz \rtimes_{-\theta} \mathbb{Z} \rightarrow 
C(S^1) \rtimes_{\theta} \mathbb{Z},\quad
z_- \mapsto U_+,\quad u_- \mapsto Z_+.\nonumber
\end{eqnarray} 
There is a natural $U(1)$-action on $C(S^3_{pq \theta})$
corresponding classically to the action in the Hopf
fibration, see \cite{hms06}. Its
restriction to the two crossed products is 
not the canonical action of $U(1)$ viewed 
as the Pontryagin dual of
$\mathbb{Z}$. However, to obtain the canonical actions
one can identify $C(S^3_{pq \theta})$ with a 
fibre product of the same crossed products, but formed 
with respect to the surjections
\begin{eqnarray}
&& \hat\pi^1_2 : \tplz \rtimes_{\theta} \mathbb{Z} \rightarrow 
C(S^1) \rtimes_{\theta} \mathbb{Z},\quad 
z_+ \mapsto Z_+,\quad u_+ \mapsto U_+,\nonumber\\ 
&& \hat\pi^2_1 : \tplz \rtimes_{-\theta} \mathbb{Z} \rightarrow 
C(S^1) \rtimes_{\theta} \mathbb{Z},\quad
z_- \mapsto Z^{-1}_+,\quad u_- \mapsto Z_+U_+.\nonumber
\end{eqnarray} 
The identification is given by
$$
	\mbox{\xymatrix@=5mm@R=16mm{
\tplz\rtimes_{\theta} \mathbb{Z} \ar[rd]+<-3mm,3mm>_{\pi^1_2} \ar@/^2pc/[rrrrr]^{\phi_1}&& 
\tplz\rtimes_{-\theta} \mathbb{Z}
\ar[ld]+<3mm,3mm>^{\pi^2_1} \ar@/^2pc/[rrrrr]^{\phi_2}&&&
\tplz\rtimes_{\theta} \mathbb{Z} \ar[rd]+<-3mm,3mm>_{\hat\pi^1_2}&& 
\tplz\rtimes_{-\theta} \mathbb{Z}\ar[ld]+<3mm,3mm>^{\hat\pi^2_1}\\
&\save[]*\txt{$C(S^1) \rtimes_{\theta} \mathbb{Z}$\mbox{\ }}
\ar[rrrrr]+<-12mm,0mm>^-{\phi_{12}}\restore & & & & &
\save[]*\txt{$C(S^1) \rtimes_{\theta} \mathbb{Z}$}\restore &
.
}}
$$
Here isomorphisms $\phi$ are given on respective generators by
$$
z\mapsto zu,\quad u\mapsto u.
$$

The $C^*$-subalgebra of $U(1)$-invariants is the $C^*$-algebra 
of the mirror quantum
2-sphere from \cite{hms06m}.
As mentioned in the introduction, we can pass from 
$C(S^3_{pq \theta})$ to the associated principal
extension, and this procedure commutes with taking
fibre products. In this way, we obtain a subalgebra 
$P \subset C(S^3_{pq \theta})$ which is a piecewise trivial 
$ \mathbb{C} \mathbb{Z}$-comodule algebra, so that it fits the setting of
this paper. The invariant subalgebra 
$P^\co$ is again the $C^*$-algebra of the mirror quantum
2-sphere.

On the other hand, there is a second natural Hopf-like $U(1)$-action 
on $C(S^3_{pq\theta})$ described in \cite{hms06m} (see also 
\cite{bhms+}). Again,
its restriction to the two crossed products making up the
fibre product $C(S^3_{pq\theta})$ is not the canonical action
of $U(1)$. This fibre product can be transformed into an isomorphic
one (carrying the canonical $U(1)$-action)
constructed by gluing two copies of $\tplz \rtimes_{-\theta} \mathbb{Z}$
over $C(S^1) \rtimes_{-\theta} \mathbb{Z}$ (with generators $Z_-,U_-$)
using the gluing maps

\begin{eqnarray}
&& \check\pi^1_2 : \tplz \rtimes_{-\theta} \mathbb{Z} \rightarrow 
C(S^1) \rtimes_{-\theta} \mathbb{Z},\quad 
z_- \mapsto Z_-,\quad u_- \mapsto U_-,\nonumber\\ 
&& \check\pi^2_1 : \tplz \rtimes_{-\theta} \mathbb{Z} \rightarrow 
C(S^1) \rtimes_{-\theta} \mathbb{Z},\quad
z_- \mapsto Z_-,\quad u_- \mapsto Z_-U_-.\nonumber
\end{eqnarray} 

The identifying maps are now given by
$$
	\mbox{\xymatrix@=5mm@R=16mm{
\tplz\rtimes_{\theta} \mathbb{Z} \ar[rd]+<-3mm,3mm>_{\pi^1_2} \ar@/^2pc/[rrrrr]^{\tilde\phi_1}&& 
\tplz\rtimes_{-\theta} \mathbb{Z}
\ar[ld]+<3mm,3mm>^{\pi^2_1} \ar@/^2pc/[rrrrr]^{\tilde\phi_2}&&&
\tplz\rtimes_{-\theta} \mathbb{Z} \ar[rd]+<-3mm,3mm>_{\check\pi^1_2}&& 
\tplz\rtimes_{-\theta} \mathbb{Z}\ar[ld]+<3mm,3mm>^{\check\pi^2_1}\\
&\save[]*\txt{$C(S^1) \rtimes_{\theta} \mathbb{Z}$\mbox{\ }}
\ar[rrrrr]+<-13mm,0mm>^-{\tilde\phi_{12}}\restore & & & & &
\save[]*\txt{$C(S^1) \rtimes_{-\theta} \mathbb{Z}$}\restore &
.
}}
$$

Here isomorphisms $\phi$ are given on generators by
\begin{eqnarray}
\tilde\phi_1&:&z_+\mapsto z_-u_-,\quad u_+\mapsto u_-^{-1},\nonumber\\
\tilde\phi_2&:&z_-\mapsto u^{-1}_-z_-,\quad u_-\mapsto u_-,\nonumber\\
\tilde\phi_{12}&:&Z_+\mapsto Z_-U_-,\quad U_+\mapsto U_-^{-1}.\nonumber
\end{eqnarray}

The subalgebra of $U(1)$-invariants is now
the
$C^*$-algebra of the 
generic Podle\'s quantum 2-sphere from \cite{p-p87}.
However, note that it is not possible to obtain the 
algebraic Podle\'s sphere in this way 
by replacing $\tplz=P_i^\co$ by the 
coordinate algebra of a 
quantum disc with generator $x$ satisfying 
$x^*x-q xx^*=1-q$ \cite{cm00}. This is related to the 
fact that already in the commutative case the algebra of
polynomial functions on a sphere has no covering corresponding
to two hemispheres -- there are no nontrivial polynomials
vanishing on a hemisphere. 
Therefore to be in this setting of fibre products
 we use more complete algebras,
e.g., $C^*$-algebras.

\section*{Acknowledgements}

The authors are very
grateful to Pawe\l\ Witkowski for producing the join picture, and to
Gabriella B\"ohm, Tomasz Brzezi\'nski, George Janelidze,
Dorota Marciniak, Tomasz Maszczyk, Fred Van Oystaeyen, Marcin 
Szamotulski and Elmar Wagner for discussions.
This work was partially supported by the
European Commission grants 
EIF-515144(UK), MKTD-CT-2004-509794 (RM), 
  the Polish Government  grants\linebreak 1~P03A~036~26 (PMH, RM),
115/E-343/SPB/6.PR UE/DIE 50/2005-2008 (PMH, BZ), and the
University of \L{}\'od\'z Grant 795 (BZ).
PMH is also grateful to the Max Planck Institute for
Mathematics in the Sciences in Leipzig for its hospitality and
financial support.

\end{document}